\documentclass[a4paper,12pt]{article}
\usepackage{amsfonts,amssymb,graphicx}

\newtheorem{lemma}{Lemma}

\begin{document}

\author{Octavian G. Mustafa\\
\textit{Faculty of Mathematics, D.A.L.,}\\
\textit{University of Craiova, Romania}\\
\textit{e-mail: octaviangenghiz@yahoo.com}\\
}
\title{Oscillatory solutions of some perturbed second order differential equations}
\date{}
\maketitle

\noindent\textbf{Abstract} We discuss the occurrence of oscillatory solutions which decay to $0$ as $s\rightarrow+\infty$ for a class of perturbed second order ordinary differential equations. As opposed to other results in the recent literature, the perturbation is as small as desired in terms of its improper integrals and it is independent of the coefficients of the non-oscillatory unperturbed equation. This class of equations reveals thus a new pathology in the theory of perturbed oscillations. 



\section{Introduction}

We are concerned in this note with a class of perturbed differential equations of second order
\begin{eqnarray}
h^{\prime\prime}+p(s,h,h^{\prime})=q(s),\qquad s\geq s_{0}>0,\label{intro_1}
\end{eqnarray}
and the occurrence of oscillatory solutions to (\ref{intro_1}) which decay to $0$ as $s\rightarrow+\infty$ when the perturbation $q$ oscillates.

Though the oscillation literature is vast, most of the investigations are generalizations of the perspective given by the fundamental papers due to A. Kartsatos \cite{Kartsatos1,Kartsatos2}, that is, if the unperturbed equation
\begin{eqnarray}
h^{\prime\prime}+p(s,h,h^{\prime})=0,\qquad s\geq s_{0},\label{intro_2}
\end{eqnarray}
is oscillatory and the function $Q$, where $Q^{\prime\prime}=q$, is itself oscillatory and decaying to $0$ when $s\rightarrow+\infty$ then the perturbed equation (\ref{intro_1}) is oscillatory. For recent advancements in this respect, the reader can consult \cite{OuWong}.

The case where the equation (\ref{intro_2}) is non-oscillatory is much more complicated, however, when the function $Q$ is strongly oscillatory in some sense, see e.g. \cite{MustafaRogovchenko}, the perturbed equation (\ref{intro_1}) will possess both oscillatory solutions, vanishing at infinity, and non-oscillatory solutions. In terms of the improper integrals used for transforming the differential equations into integral equations and in the averaging process of the former ones, this strong oscillation of $Q$ means that the integrals involved have extra large sizes in comparison with the similar integrals of the coefficients of the unperturbed equation (\ref{intro_2}).

The following problem is, to the best of our knowledge, still open: which are (if any) the conditions to be imposed on an oscillatory perturbation $q$ that will produce oscillatory, vanishing at infinity, solutions to the equation (\ref{intro_1}) whose unperturbed part (\ref{intro_2}) is not only non-oscillatory but its coefficients are also independent of $q$ ?

In this note, we present a class of second-order differential equations which are non-oscillatory but will produce such oscillatory solutions when perturbed by a small, oscillatory perturbation whose size is independent of the sizes of the coefficients of the equations.

\section{Induced oscillations to equation (\ref{intro_1})}

Consider the following perturbed equation
\begin{eqnarray}
h^{\prime\prime}+p(s)\left(h^{\prime}-\frac{h}{s}\right)+\frac{q(s)}{s}=0,\qquad s\geq s_{0}>0,\label{class_0}
\end{eqnarray}
where the coefficients $p,q:[s_{0},+\infty)\rightarrow\mathbb{R}$ are continuous and such that
\begin{eqnarray}
p(s)\geq0,\thinspace s\geq s_{0},\qquad p,q\in L^{1}((s_{0},+\infty),\mathbb{R}).\label{gen_cond_1}
\end{eqnarray}

The equation has been employed recently in studying a class of reaction-diffusion equations by means of the comparison method, see \cite{EhrnstromMustafa} and its references. It is also connected with the Lie theory of integration, see \cite[Example 2.62, p. 149]{Olver}.

According to the classical asymptotic integration theory \cite{Bellman}, the conditions (\ref{gen_cond_1}) imply that, given $c\neq0$, both the equation (\ref{class_0}) and its unperturbed part have solutions which behave asymptotically as
\begin{eqnarray}
h(s)=c\cdot s+o(s)\qquad\mbox{when }s\rightarrow+\infty,\label{asymp_integr_0}
\end{eqnarray}
so, obviously, they are non-oscillatory.

In the next sections additional restrictions will be imposed on $q$ that are independent of $p$ and will produce an oscillatory solution to the equation (\ref{class_0}) which vanishes at infinity. 

It is useful to remark that the only place with respect to the integration formula (\ref{asymp_integr_0}) where oscillations can ``hide'' is given by $c=0$. The hypotheses (\ref{gen_cond_1}), however, yield also the existence of non-oscillatory solutions to (\ref{class_0}) that obey the description (\ref{asymp_integr_0}) for $c=0$. See in this respect \cite[Proposition 1]{Mustafa2008}.

\section{Oscillation lemmas}

Consider the following ordinary differential equation
\begin{eqnarray}
\frac{dz}{ds}+p(s)z+q(s)=0,\qquad s\geq s_{0}>0,\label{ode_1}
\end{eqnarray}
where the coefficients $p,q:[s_0,+\infty)\rightarrow\mathbb{R}$ are continuous and
\begin{eqnarray*}
p(s)\geq0,\qquad s\geq s_0.
\end{eqnarray*}

\begin{lemma}\label{lem}
Assume that $p,q\in L^{1}((s_{0},+\infty),\mathbb{R})$ and there exists the unbounded from above, increasing sequence $(a_{m})_{m\geq1}$ of numbers from $[s_{0},+\infty)$ such that
\begin{eqnarray*}
q(s)<0,\thinspace s\in(a_{2m},a_{2m+1}),\qquad q(s)>0,\thinspace s\in(a_{2m+1},a_{2m+2}),
\end{eqnarray*}
and
\begin{eqnarray}
\int_{a_{2m}}^{a_{2m+1}}\vert q(s)\vert ds>3\int_{a_{2m+1}}^{a_{2m+2}}q(s)ds,\label{lem_restr_I}
\end{eqnarray}
and
\begin{eqnarray}
\int_{a_{2m+1}}^{a_{2m+2}}q(s)ds>2\int_{a_{2m+2}}^{+\infty}\vert q(s)\vert ds.\label{lem_restr_II}
\end{eqnarray}

Then, the equation (\ref{ode_1}) has an oscillatory solution $z$ such that $\lim\limits_{s\rightarrow+\infty}z(s)=0$ and
\begin{eqnarray*}
z(a_{2m})<0,\qquad z(a_{2m+1})>0
\end{eqnarray*}
for all the $m$'s great enough.
\end{lemma}

{\bf Proof.} The function $z:[s_{0},+\infty)\rightarrow\mathbb{R}$ given by the formula
\begin{eqnarray}
z(s)=\exp\left(-\int_{s_0}^{s}p(\tau)d\tau\right)\cdot\int_{s}^{+\infty}q(\tau)\exp\left(\int_{s_0}^{\tau}p(\xi)d\xi\right)d\tau,\label{form_sol}
\end{eqnarray}
where $s\geq s_0$, is a solution of equation (\ref{ode_1}).

Take $m_{0}\geq1$ great enough to verify
\begin{eqnarray}
\ln\frac{3}{2}>\int_{a_{2m}}^{+\infty}p(s)ds,\qquad m\geq m_{0}.\label{def_m_o}
\end{eqnarray}
In particular, we have the estimates
\begin{eqnarray*}
3>2\exp\left(\int_{a_{2m}}^{a_{2m+2}}p(s)ds\right)
\end{eqnarray*}
which implies that
\begin{eqnarray}
3\exp\left(\int_{s_0}^{a_{2m}}p(s)ds\right)>2\exp\left(\int_{s_0}^{a_{2m+2}}p(s)ds\right),\label{ode_1_estim_1}
\end{eqnarray}
and respectively
\begin{eqnarray*}
2>\exp\left(\int_{a_{2m+1}}^{+\infty}p(s)ds\right)
\end{eqnarray*}
which implies that
\begin{eqnarray}
2\exp\left(\int_{s_0}^{a_{2m+1}}p(s)ds\right)>\exp\left(\int_{s_0}^{+\infty}p(s)ds\right).\label{ode_1_estim_2}
\end{eqnarray}

Now, by means of (\ref{lem_restr_I}), (\ref{ode_1_estim_1}), we deduce that
\begin{eqnarray*}
&&\int_{a_{2m}}^{a_{2m+1}}\vert q(s)\vert\exp\left(\int_{s_0}^{s}p(\tau)d\tau\right)ds\\
&&\geq\int_{a_{2m}}^{a_{2m+1}}\vert q(s)\vert ds\cdot\exp\left(\int_{s_0}^{a_{2m}}p(\tau)d\tau\right)\\
&&\geq3\int_{a_{2m+1}}^{a_{2m+2}}q(s)ds\cdot\exp\left(\int_{s_0}^{a_{2m}}p(\tau)d\tau\right)\\
&&>2\exp\left(\int_{s_0}^{a_{2m+2}}p(\tau)d\tau\right)\cdot\int_{a_{2m+1}}^{a_{2m+2}}q(s)ds\\
&&\geq2\int_{a_{2m+1}}^{a_{2m+2}}q(s)\exp\left(\int_{s_0}^{s}p(\tau)d\tau\right)ds
\end{eqnarray*}
and, by means of (\ref{lem_restr_II}), (\ref{ode_1_estim_2}), that
\begin{eqnarray*}
&&\int_{a_{2m+1}}^{a_{2m+2}}q(s)\exp\left(\int_{s_0}^{s}p(\tau)d\tau\right)ds\\
&&\geq\int_{a_{2m+1}}^{a_{2m+2}}q(s)ds\cdot\exp\left(\int_{s_0}^{a_{2m+1}}p(\tau)d\tau\right)\\
&&>2\int_{a_{2m+2}}^{+\infty}\vert q(s)\vert ds\cdot\exp\left(\int_{s_0}^{a_{2m+1}}p(\tau)d\tau\right)\\
&&\geq\exp\left(\int_{s_0}^{+\infty}p(\tau)d\tau\right)\cdot\int_{a_{2m+2}}^{+\infty}\vert q(s)\vert ds\\
&&\geq\int_{a_{2m+2}}^{+\infty}\vert q(s)\vert\exp\left(\int_{s_0}^{s}p(\tau)d\tau\right) ds.
\end{eqnarray*}

Finally,
\begin{eqnarray*}
&&z(a_{2m})\cdot\exp\left(\int_{s_0}^{a_{2m}}p(\tau)d\tau\right)=\int_{a_{2m}}^{+\infty}q(\tau)\exp\left(\int_{s_0}^{\tau}p(\xi)d\xi\right)d\tau\\
&&=\left(\int_{a_{2m}}^{a_{2m+1}}+\int_{a_{2m+1}}^{a_{2m+2}}+\int_{a_{2m+2}}^{+\infty}\right)q(\tau)\exp\left(\int_{s_0}^{\tau}p(\xi)d\xi\right)d\tau\\
&&\leq\left(-\int_{a_{2m}}^{a_{2m+1}}\vert q(\tau)\vert+\int_{a_{2m+1}}^{a_{2m+2}}q(\tau)+\int_{a_{2m+2}}^{+\infty}\vert q(\tau)\vert\right)\\
&&\times\exp\left(\int_{s_0}^{\tau}p(\xi)d\xi\right)d\tau\\
&&<\left(-\int_{a_{2m}}^{a_{2m+1}}\vert q(\tau)\vert+2\int_{a_{2m+1}}^{a_{2m+2}}q(\tau)\right)\exp\left(\int_{s_0}^{\tau}p(\xi)d\xi\right)d\tau\\
&&<0
\end{eqnarray*}
and respectively
\begin{eqnarray*}
&&z(a_{2m+1})\exp\left(\int_{s_0}^{a_{2m+1}}p(\tau)d\tau\right)=\int_{a_{2m+1}}^{+\infty}q(\tau)\exp\left(\int_{s_0}^{\tau}p(\xi)d\xi\right)d\tau\\
&&=\left(\int_{a_{2m+1}}^{a_{2m+2}}+\int_{a_{2m+2}}^{+\infty}\right)q(\tau)\exp\left(\int_{s_0}^{\tau}p(\xi)d\xi\right)d\tau\\
&&\geq\left(\int_{a_{2m+1}}^{a_{2m+2}}q(\tau)-\int_{a_{2m+2}}^{+\infty}\vert q(\tau)\vert\right)\exp\left(\int_{s_0}^{\tau}p(\xi)d\xi\right)d\tau\\
&&>0.
\end{eqnarray*}

The proof is complete. $\square$

\begin{lemma}\label{lem_2}
Assume that the hypotheses of Lemma \ref{lem} hold. Suppose also that
\begin{eqnarray}
\int_{a_{2m}}^{a_{2m+1}}(s-a_{2m})\vert q(s)\vert ds>6\cdot\frac{a_{2m+1}^{2}}{a_{2m}}\int_{a_{2m+1}}^{+\infty}\vert q(s)\vert ds\label{lem_2_restr_1}
\end{eqnarray}
and that
\begin{eqnarray}
\int_{a_{2m+1}}^{a_{2m+2}}(s-a_{2m+1})q(s) ds>4\cdot\frac{a_{2m+2}^{2}}{a_{2m+1}}\int_{a_{2m+2}}^{+\infty}\vert q(s)\vert ds\label{lem_2_restr_2}
\end{eqnarray}
for all the $m$'s great enough.

Then, the function $h:[s_{0},+\infty)\rightarrow\mathbb{R}$ with the formula
\begin{eqnarray}
h(s)=-s\int_{s}^{+\infty}\frac{z(\tau)}{\tau^2}d\tau,\qquad s\geq s_{0},\label{lem_2_form_h}
\end{eqnarray}
where $z$ is the oscillatory solution (\ref{form_sol}), is itself oscillatory, that is
\begin{eqnarray*}
h(a_{2m})>0,\qquad h(a_{2m+1})<0
\end{eqnarray*}
for all the $m$'s great enough.
\end{lemma}

{\bf Proof.} Recall the definition (\ref{def_m_o}) of $m_0$.

{\it Step 1.} We shall establish that $h(a_{2m})>0$.

As before,
\begin{eqnarray}
\int_{a_{2m}}^{+\infty}\frac{z(s)}{s^2}ds=\left(\int_{a_{2m}}^{a_{2m+1}}+\int_{a_{2m+1}}^{a_{2m+2}}+\int_{a_{2m+2}}^{+\infty}\right)\frac{z(s)}{s^2}ds.\label{dec_1}
\end{eqnarray}

Notice that the third term of the decomposition can be estimated by
\begin{eqnarray*}
\left\vert\int_{a_{2m+2}}^{+\infty}\frac{z(s)}{s^2}ds\right\vert&\leq&\int_{a_{2m+2}}^{+\infty}\frac{\vert z(s)\vert}{s^2}ds\\
&\leq&\int_{a_{2m+2}}^{+\infty}\frac{1}{s^2}\exp\left(\int_{s}^{+\infty}p(\tau)d\tau\right)\int_{s}^{+\infty}\vert q(\tau)\vert d\tau ds\\
&\leq&\int_{a_{2m+2}}^{+\infty}\frac{ds}{s^2}\cdot\exp\left(\int_{a_{2m+2}}^{+\infty}p(\tau)d\tau\right)\cdot\int_{a_{2m+2}}^{+\infty}\vert q(\tau)\vert d\tau\\
&<&\frac{3/2}{a_{2m+2}}\int_{a_{2m+2}}^{+\infty}\vert q(\tau)\vert d\tau<\frac{2}{a_{2m+1}}\int_{a_{2m+2}}^{+\infty}\vert q(\tau)\vert d\tau,
\end{eqnarray*}
where $m\geq m_{0}$.

The first term of the decomposition (\ref{dec_1}) can be decomposed further into
\begin{eqnarray*}
\int_{a_{2m}}^{a_{2m+1}}\frac{z(s)}{s^2}ds&=&\int_{a_{2m}}^{a_{2m+1}}\frac{1}{s^2}\exp\left(-\int_{s_0}^{s}p(\tau)d\tau\right)\\
&\times&\left(\int_{s}^{a_{2m+1}}+\int_{a_{2m+1}}^{+\infty}\right)q(\tau)\exp\left(\int_{s_0}^{\tau}p(\xi)d\xi\right)d\tau ds.
\end{eqnarray*}

The second part of this decomposition is estimated by
\begin{eqnarray*}
&&\left\vert\int_{a_{2m}}^{a_{2m+1}}\frac{1}{s^2}\exp\left(-\int_{s_0}^{s}p(\tau)d\tau\right)\int_{a_{2m+1}}^{+\infty}q(\tau)\exp\left(\int_{s_0}^{\tau}p(\xi)d\xi\right)d\tau ds\right\vert\\
&&\leq\int_{a_{2m}}^{a_{2m+1}}\frac{ds}{s^2}\cdot\exp\left(\int_{a_{2m}}^{+\infty}p(\tau)d\tau\right)\cdot\int_{a_{2m+1}}^{+\infty}\vert q(\tau)\vert d\tau\\
&&<\frac{3/2}{a_{2m}}\int_{a_{2m+1}}^{+\infty}\vert q(\tau)\vert d\tau<\frac{2}{a_{2m}}\int_{a_{2m+1}}^{+\infty}\vert q(\tau)\vert d\tau,
\end{eqnarray*}
where $m\geq m_{0}$, while the first part is a negative quantity.

The second term of the decomposition (\ref{dec_1}) is estimated by
\begin{eqnarray*}
&&\left\vert\int_{a_{2m+1}}^{a_{2m+2}}\frac{z(s)}{s^2}ds\right\vert\leq\int_{a_{2m+1}}^{a_{2m+2}}\frac{ds}{s^2}\cdot\exp\left(\int_{a_{2m+1}}^{+\infty}p(\tau)d\tau\right)\cdot\int_{a_{2m+1}}^{+\infty}\vert q(\tau)\vert d\tau\\
&&<\frac{3/2}{a_{2m+1}}\int_{a_{2m+1}}^{+\infty}\vert q(\tau)\vert d\tau<\frac{2}{a_{2m+1}}\int_{a_{2m+1}}^{+\infty}\vert q(\tau)\vert d\tau,
\end{eqnarray*}
where $m\geq m_{0}$.

Collecting all the estimates, we deduce that
\begin{eqnarray*}
&&\int_{a_{2m}}^{+\infty}\frac{z(s)}{s^2}ds\\
&&<\int_{a_{2m}}^{a_{2m+1}}\frac{1}{s^2}\exp\left(-\int_{s_0}^{s}p(\tau)d\tau\right)\int_{s}^{a_{2m+1}}q(\tau)\exp\left(\int_{s_0}^{\tau}p(\xi)d\xi\right)d\tau ds\\
&&+\left(\frac{2}{a_{2m}}+\frac{2}{a_{2m+1}}\right)\int_{a_{2m+1}}^{+\infty}\vert q(\tau)\vert d\tau+\frac{2}{a_{2m+1}}\int_{a_{2m+2}}^{+\infty}\vert q(\tau)\vert d\tau\\
&&<\int_{a_{2m}}^{a_{2m+1}}\frac{1}{s^2}\exp\left(-\int_{s_0}^{s}p(\tau)d\tau\right)\int_{s}^{a_{2m+1}}q(\tau)\exp\left(\int_{s_0}^{\tau}p(\xi)d\xi\right)d\tau ds\\
&&+\frac{6}{a_{2m}}\int_{a_{2m+1}}^{+\infty}\vert q(\tau)\vert d\tau.
\end{eqnarray*}

Finally, notice that
\begin{eqnarray*}
&&\left\vert\int_{a_{2m}}^{a_{2m+1}}\frac{1}{s^2}\exp\left(-\int_{s_0}^{s}p(\tau)d\tau\right)\int_{s}^{a_{2m+1}}q(\tau)\exp\left(\int_{s_0}^{\tau}p(\xi)d\xi\right)d\tau ds\right\vert\\
&&=\int_{a_{2m}}^{a_{2m+1}}\frac{1}{s^2}\exp\left(-\int_{s_0}^{s}p(\tau)d\tau\right)\int_{s}^{a_{2m+1}}\vert q(\tau)\vert\exp\left(\int_{s_0}^{\tau}p(\xi)d\xi\right)d\tau ds\\
&&\geq\int_{a_{2m}}^{a_{2m+1}}\frac{1}{s^2}\exp\left(-\int_{s_0}^{s}p(\tau)d\tau\right)\int_{s}^{a_{2m+1}}\vert q(\tau)\vert d\tau\exp\left(\int_{s_0}^{s}p(\xi)d\xi\right)ds\\
&&\geq\int_{a_{2m}}^{a_{2m+1}}\frac{1}{s^2}\int_{s}^{a_{2m+1}}\vert q(\tau)\vert d\tau ds=\int_{a_{2m}}^{a_{2m+1}}\frac{d}{ds}\left(-\frac{1}{s}\right)\int_{s}^{a_{2m+1}}\vert q(\tau)\vert d\tau ds\\
&&=\frac{1}{a_{2m}}\int_{a_{2m}}^{a_{2m+1}}\vert q(\tau)\vert d\tau-\int_{a_{2m}}^{a_{2m+1}}\frac{\vert q(s)\vert}{s}ds\\
&&=\int_{a_{2m}}^{a_{2m+1}}\frac{s-a_{2m}}{s}\vert q(s)\vert ds\cdot\frac{1}{a_{2m}}\\
&&>\frac{1}{a_{2m+1}a_{2m}}\cdot\int_{a_{2m}}^{a_{2m+1}}(s-a_{2m})\vert q(s)\vert ds\\
&&>\frac{1}{a_{2m+1}a_{2m}}\cdot6\frac{a_{2m+1}^{2}}{a_{2m}}\int_{a_{2m+1}}^{+\infty}\vert q(s)\vert ds=6\frac{a_{2m+1}}{a_{2m}^{2}}\int_{a_{2m+1}}^{+\infty}\vert q(s)\vert ds\\
&&>\frac{6}{a_{2m}}\int_{a_{2m+1}}^{+\infty}\vert q(s)\vert ds.
\end{eqnarray*}

In conclusion, since $q$ is negative valued in $(a_{2m},a_{2m+1})$, we get that
\begin{eqnarray*}
&&\int_{a_{2m}}^{+\infty}\frac{z(s)}{s^2}ds\\
&&<-\int_{a_{2m}}^{a_{2m+1}}\frac{1}{s^2}\exp\left(-\int_{s_0}^{s}p(\tau)d\tau\right)\int_{s}^{a_{2m+1}}\vert q(\tau)\vert\exp\left(\int_{s_0}^{\tau}p(\xi)d\xi\right)d\tau ds\\
&&+\frac{6}{a_{2m}}\int_{a_{2m+1}}^{+\infty}\vert q(\tau)\vert d\tau\\
&&<0.
\end{eqnarray*}

The first step is complete.

{\it Step 2.} We shall establish that $h(a_{2m+1})<0$.

First,
\begin{eqnarray}
\int_{a_{2m+1}}^{+\infty}\frac{z(s)}{s^2}ds=\left(\int_{a_{2m+1}}^{a_{2m+2}}+\int_{a_{2m+2}}^{+\infty}\right)\frac{z(s)}{s^2}ds.\label{dec_2}
\end{eqnarray}

The second term of the decomposition is estimated by
\begin{eqnarray*}
&&\left\vert\int_{a_{2m+2}}^{+\infty}\frac{z(s)}{s^2}ds\right\vert\\
&&\leq\int_{a_{2m+2}}^{+\infty}\frac{ds}{s^2}\cdot\exp\left(\int_{a_{2m+2}}^{+\infty}p(\tau)d\tau\right)\cdot\int_{a_{2m+2}}^{+\infty}\vert q(s)\vert ds\\
&&<\frac{2}{a_{2m+2}}\int_{a_{2m+2}}^{+\infty}\vert q(s)\vert ds,
\end{eqnarray*}
where $m\geq m_{0}$.

The first term of the decomposition (\ref{dec_2}) is decomposed further into
\begin{eqnarray*}
&&\int_{a_{2m+1}}^{a_{2m+2}}\frac{z(s)}{s^2}ds=\int_{a_{2m+1}}^{a_{2m+2}}\frac{1}{s^2}\exp\left(-\int_{s_0}^{s}p(\tau)d\tau\right)\\
&&\times\left(\int_{s}^{a_{2m+2}}+\int_{a_{2m+2}}^{+\infty}\right)q(\tau)\exp\left(\int_{s_0}^{\tau}p(\xi)d\xi\right)d\tau ds.
\end{eqnarray*}

The second part of this decomposition is estimated by
\begin{eqnarray*}
&&\left\vert\int_{a_{2m+1}}^{a_{2m+2}}\frac{1}{s^2}\exp\left(-\int_{s_0}^{s}p(\tau)d\tau\right)\int_{a_{2m+2}}^{+\infty}q(\tau)\exp\left(\int_{s_0}^{\tau}p(\xi)d\xi\right)d\tau ds\right\vert\\
&&\leq\int_{a_{2m+1}}^{a_{2m+2}}\frac{ds}{s^2}\cdot\exp\left(\int_{a_{2m+1}}^{+\infty}p(\tau)d\tau\right)\cdot\int_{a_{2m+2}}^{+\infty}\vert q(s)\vert ds\\
&&<\frac{2}{a_{2m+1}}\int_{a_{2m+2}}^{+\infty}\vert q(s)\vert ds,
\end{eqnarray*}
where $m\geq m_{0}$.

Finally, notice that
\begin{eqnarray*}
&&\int_{a_{2m+1}}^{a_{2m+2}}\frac{1}{s^2}\exp\left(-\int_{s_0}^{s}p(\tau)d\tau\right)\int_{s}^{a_{2m+2}}q(\tau)\exp\left(\int_{s_0}^{\tau}p(\xi)d\xi\right)d\tau ds\\
&&\geq\int_{a_{2m+1}}^{a_{2m+2}}\frac{1}{s^2}\exp\left(-\int_{s_0}^{s}p(\tau)d\tau\right)\int_{s}^{a_{2m+2}} q(\tau) d\tau\exp\left(\int_{s_0}^{s}p(\xi)d\xi\right)ds\\
&&\geq\int_{a_{2m+1}}^{a_{2m+2}}\frac{1}{s^2}\int_{s}^{a_{2m+2}} q(\tau) d\tau ds\\
&&\geq\int_{a_{2m+1}}^{a_{2m+2}}(s-a_{2m+1}) q(s) ds\cdot\frac{1}{a_{2m+2}a_{2m+1}}\\
&&>4\frac{a_{2m+2}}{a_{2m+1}^{2}}\int_{a_{2m+2}}^{+\infty}\vert q(s)\vert ds>\frac{4}{a_{2m+1}}\int_{a_{2m+2}}^{+\infty}\vert q(s)\vert ds.
\end{eqnarray*}

In conclusion, since $q$ is positive valued in $(a_{2m+1},a_{2m+2})$, we get that
\begin{eqnarray*}
&&\int_{a_{2m+1}}^{+\infty}\frac{z(s)}{s^2}ds\\
&&>\int_{a_{2m+1}}^{a_{2m+2}}\frac{1}{s^2}\exp\left(-\int_{s_0}^{s}p(\tau)d\tau\right)\int_{s}^{a_{2m+2}}q(\tau)\exp\left(\int_{s_0}^{\tau}p(\xi)d\xi\right)d\tau ds\\
&&-\frac{4}{a_{2m+1}}\int_{a_{2m+2}}^{+\infty}\vert q(s)\vert ds\\
&&>0.
\end{eqnarray*}

The second step is complete. $\square$

Observe that
\begin{eqnarray*}
\frac{d}{ds}\left(\frac{h(s)}{s}\right)=\frac{z(s)}{s^2},\qquad s\geq s_{0}.
\end{eqnarray*}

Concluding this section, the function $h$ given by
\begin{eqnarray*}
\frac{h(s)}{s}=-\int_{s}^{+\infty}\left(\frac{h(\tau)}{\tau}\right)^{\prime}d\tau=-\int_{s}^{+\infty}\frac{z(\tau)}{\tau^{2}}d\tau
\end{eqnarray*}
is an oscillatory solution of the equation (\ref{class_0}). Since, according to Lemma \ref{lem}, $\lim\limits_{s\rightarrow+\infty}z(s)=0$, an application of the L'H\^{o}pital rule yields
\begin{eqnarray*}
\lim\limits_{s\rightarrow+\infty}h(s)=\lim\limits_{s\rightarrow+\infty}s\int_{s}^{+\infty}\frac{z(\tau)}{\tau^{2}}d\tau=0.
\end{eqnarray*}

\section{Example}

We start with an auxiliary result.

\begin{lemma}\label{lem_3}
Set $C>0$. There exists $\varepsilon\in(0,1)$ small enough such that
\begin{eqnarray}
C\cdot\frac{\varepsilon^{m^2}}{m}>\sum\limits_{k=m+1}^{+\infty}k\varepsilon^{k^2},\qquad m\geq1.\label{aux_ineq}
\end{eqnarray}
\end{lemma}

{\bf Proof.} The inequality can be recast as
\begin{eqnarray}
\frac{C}{m}>\sum\limits_{p=1}^{+\infty}(m+p)\varepsilon^{(m+p)^{2}-m^{2}},\qquad m\geq1.\label{aux_ineq_aux}
\end{eqnarray}

Notice that $(m+p)^{2}-m^{2}=p^{2}+2mp\geq2m+p$ for all $m,p\geq1$.

Since $\varepsilon\in(0,1)$, we deduce that
\begin{eqnarray*}
\sum\limits_{p=1}^{+\infty}(m+p)\varepsilon^{(m+p)^{2}-m^{2}}&<&\varepsilon^{m}\cdot\varepsilon\sum\limits_{p=1}^{+\infty}(m+p)\varepsilon^{m+p-1}\\
&<&\varepsilon^{m}\cdot\varepsilon\sum\limits_{q=1}^{+\infty}q\varepsilon^{q-1}=\varepsilon^{m}\cdot\frac{\varepsilon}{(1-\varepsilon)^2},
\end{eqnarray*}
where $m\geq1$.

Notice that the function $x\mapsto x\varepsilon^{x}$ is decreasing in $[1,+\infty)$ when
\begin{eqnarray*}
\ln\left(\frac{1}{\varepsilon}\right)>1.
\end{eqnarray*}
In particular, this holds for $\varepsilon\in\left(0,\frac{1}{3}\right)$.

By replacing the inequalities (\ref{aux_ineq}), (\ref{aux_ineq_aux}) with the sharper one below
\begin{eqnarray}
C>\varepsilon\cdot\frac{\varepsilon}{(1-\varepsilon)^{2}}\geq m\varepsilon^{m}\cdot\frac{\varepsilon}{(1-\varepsilon)^{2}},\label{aux_ineq_2}
\end{eqnarray}
we are able to estimate the size of $\varepsilon$ from the restrictions
\begin{eqnarray*}
1-2\varepsilon>\left(\frac{1}{C}-1\right)\varepsilon^{2},\qquad \varepsilon\in\left(0,\frac{1}{3}\right).
\end{eqnarray*}

Observe that the left-hand part of the inequality is greater than $\varepsilon$.

So, if $C\geq1$ then the inequality (\ref{aux_ineq_2}) holds everywhere in $\left(0,\frac{1}{3}\right)$. If $C\in(0,1)$ then, by replacing (\ref{aux_ineq_2}) with the stronger inequality
\begin{eqnarray*}
\varepsilon>\left(\frac{1}{C}-1\right)\varepsilon^{2},
\end{eqnarray*}
we conclude that the inequality (\ref{aux_ineq_2}) holds for every $\varepsilon\in(0,\varepsilon_0)$, where
\begin{eqnarray*}
\varepsilon_0=\min\left\{\frac{1}{3},\frac{C}{1-C}\right\},
\end{eqnarray*}
and so does (\ref{aux_ineq}).

The proof is complete. $\square$

Set now $\varepsilon\in\left(0,\frac{1}{3}\right)$. Further restrictions will be imposed on $\varepsilon$ in the course of our construction.

Introduce the sequences $(a_{m})_{m\geq1}$, $(c_{m}(\alpha))_{m\geq1}$, $(d_{m}(\beta))_{m\geq1}$ via the formulas
\begin{eqnarray*}
a_{m}=m\pi,\quad c_{m}(\alpha)=\alpha\cdot\varepsilon^{m^2},\quad d_{m}(\beta)=\beta\cdot\frac{\varepsilon^{m^2}}{m},\quad m\geq1,
\end{eqnarray*}
for some $\alpha,\beta>0$.

Further, define the function $q(\alpha,\beta):[s_{0},+\infty)\rightarrow\mathbb{R}$ of class $C^{1}$ such that
\begin{eqnarray}
q(\alpha,\beta)(s)=\left\{
\begin{array}{ll}
-c_{m}(\alpha)\sin^{2}s,\thinspace s\in[a_{2m},a_{2m+1}],\\
d_{m}(\beta)\sin^{2}s,\thinspace s\in[a_{2m+1},a_{2m+2}],
\end{array}
\right.\label{q_alfa_beta}
\end{eqnarray}
for all $m\geq m_{1}=\max\{1,\frac{s_{0}}{2\pi}\}$.

Observe that
\begin{eqnarray}
&&\int_{a_{2m}}^{a_{2m+1}}(s-a_{2m})\vert q(\alpha,\beta)(s)\vert ds=\int_{a_{2m}}^{a_{2m+1}}\int^{a_{2m+1}}_{s}\vert q(\alpha,\beta)(\tau)\vert d\tau ds\nonumber\\
&&=c_{m}(\alpha)\int_{a_{2m}}^{a_{2m+1}}\left(\frac{a_{2m+1}-s}{2}+\frac{\sin2s}{4}\right)ds\nonumber\\
&&=\frac{c_m(\alpha)}{4}\left[(a_{2m+1}-a_{2m})^{2}-\frac{\cos2a_{2m+1}-\cos2a_{2m}}{2}\right]\nonumber\\
&&=\frac{c_m(\alpha)}{4}(a_{2m+1}-a_{2m})^{2}=\frac{\pi^2}{4}\cdot c_{m}(\alpha)\label{ex_integr_1}
\end{eqnarray}
and respectively
\begin{eqnarray}
\int_{a_{2m+1}}^{a_{2m+2}}(s-a_{2m+1}) q(\alpha,\beta)(s)ds=\frac{\pi^2}{4}\cdot d_{m}(\beta),\label{ex_integr_2}
\end{eqnarray}
where $m\geq m_{1}$.

We have also the next estimates
\begin{eqnarray*}
&&\int_{a_{2m+1}}^{+\infty}s\vert q(\alpha,\beta)(s)\vert ds\\
&&=\sum\limits_{k=m}^{+\infty}\left(\int_{a_{2k+1}}^{a_{2k+2}}+\int_{a_{2k+2}}^{a_{2k+3}}\right)s\vert q(\alpha,\beta)(s)\vert ds\\
&&\leq\sum\limits_{k=m}^{+\infty}\left(\int_{a_{2k+1}}^{a_{2k+2}}a_{2k+2}d_{k}(\beta)ds+\int_{a_{2k+2}}^{a_{2k+3}}a_{2k+3}c_{k+1}(\alpha)ds\right)\\
&&=\pi^{2}\sum\limits_{k=m}^{+\infty}[(2k+2)d_{k}(\beta)+(2k+3)c_{k+1}(\alpha)]\\
&&=\pi^{2}\left\{(2m+2)d_{m}(\beta)+\sum\limits_{k=m+1}^{+\infty}[(2k+2)d_{k}(\beta)+(2k+1)c_{k}(\alpha)]\right\}.
\end{eqnarray*}

Assume that $\alpha>3\beta$. Remark the obvious relations
\begin{eqnarray*}
c_{k}(\alpha)=k\frac{\alpha}{\beta}\cdot d_{k}(\beta)>3d_{k}(\beta),\qquad k\geq1.
\end{eqnarray*}

So,
\begin{eqnarray*}
&&\int_{a_{2m+1}}^{+\infty}s\vert q(\alpha,\beta)(s)\vert ds\\
&&\leq\pi^{2}\left\{4\frac{\beta}{\alpha}c_{m}(\alpha)+\sum\limits_{k=m+1}^{+\infty}\left[4\frac{\beta}{\alpha}c_{k}(\alpha)+(2k+1)c_{k}(\alpha)\right]\right\}\\
&&<\pi^{2}\left[4\frac{\beta}{\alpha}c_{m}(\alpha)+\sum\limits_{k=m+1}^{+\infty}(3k+3)c_{k}(\alpha)\right]\\
&&\leq\pi^{2}\left(4\frac{\beta}{\alpha}c_{m}(\alpha)+6\sum\limits_{k=m+1}^{+\infty}kc_{k}(\alpha)\right).
\end{eqnarray*}

Let us diminish the ratio $\frac{\beta}{\alpha}$ by asking that 
\begin{eqnarray}
\alpha>192\beta.\label{ratio_alfa_beta}
\end{eqnarray}

Thus, the next inequality
\begin{eqnarray*}
\frac{1}{6}\left(\frac{1}{48}-\frac{4\beta}{\alpha}\right)\cdot c_{m}(\alpha)>\sum\limits_{k=m+1}^{+\infty}kc_{k}(\alpha),\quad m\geq1,
\end{eqnarray*}
follows from Lemma \ref{lem_3} provided that $\varepsilon<\varepsilon_{0}=\varepsilon_{0}(C)=\varepsilon_{0}(\alpha,\beta)$.

The meaning of this inequality is given by the estimates
\begin{eqnarray*}
&&\int_{a_{2m}}^{a_{2m+1}}(s-a_{2m})\vert q(\alpha,\beta)(s)\vert ds=\int_{a_{2m}}^{a_{2m+1}}\int^{a_{2m}}_{s}\vert q(\alpha,\beta)(\tau)\vert d\tau ds\\
&&=\frac{\pi^2}{4}\cdot c_{m}(\alpha)\\
&&>12\pi^{2}\left(4\frac{\beta}{\alpha}c_{m}(\alpha)+6\sum\limits_{k=m+1}^{+\infty}kc_{k}(\alpha)\right)\\
&&>12\int_{a_{2m+1}}^{+\infty}s\vert q(\alpha,\beta)(s)\vert ds>6\cdot\frac{a_{2m+1}}{a_{2m}}\int_{a_{2m+1}}^{+\infty}s\vert q(\alpha,\beta)(s)\vert ds\\
&&>6\cdot\frac{a_{2m+1}^{2}}{a_{2m}}\int_{a_{2m+1}}^{+\infty}\vert q(\alpha,\beta)(s)\vert ds,
\end{eqnarray*}
that is, our example satisfies the restriction (\ref{lem_2_restr_1}).

The size of $\varepsilon$ will be reduced again. In fact, consider the inequality
\begin{eqnarray*}
\frac{1}{192}\cdot\frac{\beta}{\alpha}\cdot\frac{c_{m}(\alpha)}{m}>\sum\limits_{k=m+1}^{+\infty}kc_{k}(\alpha),\quad m\geq1.
\end{eqnarray*}
The magnitude of $\varepsilon$ follows now from Lemma \ref{lem_3}.

The meaning of the inequality is given by the estimates
\begin{eqnarray*}
&&\int_{a_{2m+1}}^{a_{2m+2}}(s-a_{2m+1}) q(\alpha,\beta)(s)ds=\frac{\pi^2}{4}\cdot d_{m}(\beta)\\
&&>48\pi^{2}\sum\limits_{k=m+1}^{+\infty}kc_{k}(\alpha)\\
&&>8\int_{a_{2m+2}}^{+\infty}s\vert q(\alpha,\beta)(s)\vert ds\\
&&>4\cdot\frac{a_{2m+2}}{a_{2m+1}}\int_{a_{2m+2}}^{+\infty}s\vert q(\alpha,\beta)(s)\vert ds\\
&&>4\cdot\frac{a_{2m+2}^{2}}{a_{2m+1}}\int_{a_{2m+2}}^{+\infty}\vert q(\alpha,\beta)(s)\vert ds,
\end{eqnarray*}
so, our example satisfies the restriction (\ref{lem_2_restr_2}).

Recalling the computations (\ref{ex_integr_1}), (\ref{ex_integr_2}), we remark that
\begin{eqnarray*}
\int_{a_{2m}}^{a_{2m+1}}\vert q(\alpha,\beta)(s)\vert ds&=&\frac{\pi}{2}\cdot c_{m}(\alpha)\\
&=&\frac{2}{\pi}\int_{a_{2m}}^{a_{2m+1}}(s-a_{2m})\vert q(\alpha,\beta)(s)\vert ds\\
&>&\frac{2}{\pi}\cdot12\int_{a_{2m+1}}^{+\infty}s\vert q(\alpha,\beta)(s)\vert ds\\
&>&6\int_{a_{2m+1}}^{+\infty}s\vert q(\alpha,\beta)(s)\vert ds\\
&>&3\int_{a_{2m+1}}^{a_{2m+2}} q(\alpha,\beta)(s) ds,
\end{eqnarray*}
which is the same as (\ref{lem_restr_I}).

Further, we have
\begin{eqnarray*}
\int_{a_{2m+1}}^{a_{2m+2}} q(\alpha,\beta)(s) ds&=&\frac{\pi}{2}\cdot d_{m}(\beta)\\
&=&\frac{2}{\pi}\int_{a_{2m+1}}^{a_{2m+2}}(s-a_{2m+1}) q(\alpha,\beta)(s) ds\\
&>&\frac{2}{\pi}\cdot8\int_{a_{2m+2}}^{+\infty}s\vert q(\alpha,\beta)(s)\vert ds\\
&>&4\int_{a_{2m+2}}^{+\infty}s\vert q(\alpha,\beta)(s)\vert ds\\
&>&2\int_{a_{2m+2}}^{+\infty}\vert q(\alpha,\beta)(s)\vert ds,
\end{eqnarray*}
which is the same as (\ref{lem_restr_II}).


\begin{thebibliography}{50}

\bibitem{Bellman} R. Bellman, Stability theory of differential equations, McGraw-Hill, New York, 1953.

\bibitem{EhrnstromMustafa} M. Ehrnstr\"{o}m, O.G. Mustafa, On positive solutions of a class of nonlinear elliptic equations, Nonlinear Anal. TMA 67 (2007), 1147--1154.

\bibitem{Kartsatos1} A.G. Kartsatos, On the maintenance of oscillations of $n$-th order equations under a small forcing term, J. Differential Equations 10 (1971), 355--363.

\bibitem{Kartsatos2} A.G. Kartsatos, The oscillation of a forced equation implies the oscillation of the unforced equation --- small forcing, J. Math. Anal. Appl. 76 (1980), 98--106.

\bibitem{MustafaRogovchenko} O.G. Mustafa, Yu.V. Rogovchenko, Oscillation of second-order perturbed differential equations, Math. Nachr. 278 (2005), 460--469.

\bibitem{Mustafa2008} O.G. Mustafa, Existence of positive evanescent solutions to some quasilinear elliptic equations, Bull. Austral. Math. Soc. 78 (2008), 157--162.

\bibitem{Olver} P.J. Olver, Applications of Lie groups to differential equations. Second edition, Springer-Verlag, New York, 2000.

\bibitem{OuWong} C.H. Ou, J.S.W. Wong, Forced oscillation of $n$-th order functional differential equations, J. Math. Anal. Appl. 262 (2001), 722--732.



\end{thebibliography}
\end{document}